\documentclass[11pt,a4paper]{article}

\usepackage{amsmath,amssymb,amsthm}
\usepackage{mathtools}
\usepackage{enumitem}
\usepackage[hidelinks]{hyperref}
\usepackage{microtype}
\usepackage{authblk}

\newtheorem{theorem}{Theorem}[section]
\newtheorem{proposition}[theorem]{Proposition}
\newtheorem{corollary}[theorem]{Corollary}
\newtheorem{lemma}[theorem]{Lemma}
\theoremstyle{remark}
\newtheorem{remark}[theorem]{Remark}

\newcommand{\cl}{\operatorname{cl}}

\newcommand{\lsupl}{\operatorname*{limsup}_{\ell}}
\newcommand{\lsupo}{\operatorname*{limsup}_{o}}

\title{Sequential convergence in the lower topology and co-sobriety of Scott spaces}

\author[1]{Xiaoquan Xu\thanks{Corresponding author: \texttt{xiqxu2002@163.com}}}

\author[2]{Wei Ji\thanks{Email: \texttt{javeey@163.com}}}

\affil[1]{School of Computer Information Engineering, Nanchang Institute of Technology, Nanchang 330044, China}
\affil[2]{School of Mathematics and Statistics, Guilin University of Technology, Guilin 541004, China}
\date{}

\begin{document}
\maketitle

\begin{abstract}
We investigate how order structure controls sequential convergence in the lower topology of a poset and how this, in turn, governs topological duality and co-sobriety of Scott spaces. We first give a purely order-theoretic characterization of lower convergence. For complete lattices this yields a lower-limsup formula, and on complete chains the lower limsup coincides with the usual lattice limsup. We also characterize first countability of the lower topology at a point by the existence of a countable coinitial subset outside its principal ideal. We prove that every compact saturated subset of a Scott space is sequentially closed in the lower topology and deduce that, for every poset, the lower topology and the co-compact (de Groot dual) topology of its Scott space have exactly the same convergent sequences and limits. For coherent Scott spaces we identify the lower closure of a $k$-irreducible compact saturated set with its double-cone closure and obtain an infimum criterion for co-sobriety. Consequently, every bounded-complete dcpo with sequential lower topology has a co-sober Scott space. In particular, every countable bounded-complete dcpo, and hence every countable complete lattice, has a co-sober Scott space.
\end{abstract}



\smallskip
\noindent\textbf{Funding:} This research was supported by the National Natural Science Foundation of China (Nos. 12471070, 12071199, 12561083).

\section{Introduction}

Co-sobriety was introduced by Escard\'o, Lawson and Simpson in their study of a dual form of the Hofmann--Mislove theorem \cite{Escardo-Lawson-Simpson}. A $T_0$-space is co-sober if every $k$-irreducible compact saturated subset is the saturation of a point. Co-sobriety and sobriety are distinct: sober spaces need not be co-sober, and even a sober Scott space can fail to be co-sober \cite{Shen-Wu-Xi-Zhao,Wen-Xu}.

Let $P$ be a poset. We write $\omega(P)$ for its lower topology, $\sigma(P)$ for its Scott topology, and $\Sigma P=(P,\sigma(P))$ for its Scott space. The co-compact, or de Groot dual, topology of $\Sigma P$ will be denoted by $\sigma(P)^d$; it is generated by the complements of the compact saturated subsets of $\Sigma P$.

Topology-level duality between one-sided order topologies and de Groot duals has been studied under stable compactness and related structural hypotheses \cite{Lawson-stably-compact,Xu-duals,Xu-order-topology}. In particular, systematic criteria for when the de Groot dual of a Scott topology is the lower topology are developed in \cite[Chapters 9--10]{Xu-order-topology}. Lawson and Xu developed a systematic account of the lower topology on general $T_0$-spaces and its interaction with compactness, well-filteredness and coherence \cite{Lawson-Xu-lower}. Zhao, He and Wang proved, by a different method, that every locally compact coherent $T_0$-space is co-sober \cite{Zhao-He-Wang}. The present paper takes a different viewpoint, centered on the underlying order. Rather than assuming global topological hypotheses that force topology-level duality, we begin with the order relation itself and use ordinary sequences as the organizing principle. This leads from an order-theoretic description of lower convergence to topological consequences for de Groot duality and co-sobriety.

Our first result makes the order-theoretic nature of lower convergence explicit. For a poset $P$ and a sequence $(x_n)_{n\geq 1}$ in $P$, $(x_n)_{n\geq 1}$ converges to $x$ in $\omega(P)$ precisely when every element not below $x$ is eventually not below $x_n$. Equivalently,
\[
  \bigcap_{m\ge 1}{\downarrow}\{x_n:n\ge m\}\subseteq{\downarrow} x.
\]
For a complete lattice this gives a natural lower limsup. If
\[
 \lsupl x_n
 =\bigvee\bigcap_{m\ge1}{\downarrow}\{x_n:n\ge m\},
\]
then
\[
 x_n\longrightarrow x\text{ in }\omega(L)
 ~\Longleftrightarrow~
 \lsupl x_n\le x.
\]
On a complete chain, $\lsupl x_n$ agrees with the usual lattice limsup
\[
 \bigwedge_{m\ge1}\bigvee_{n\ge m}x_n.
\]
Thus, for the extended real line with its usual order, lower-topological convergence is exactly the one-sided condition $\lsupl x_n\le x$.

The same order-theoretic viewpoint also identifies the precise local countability condition behind this convergence. The lower topology is first-countable at $x$ if and only if $P\setminus{\downarrow} x$ has a countable coinitial subset. This supplies a purely order-theoretic source of first countability, and hence of the Fr\'echet--Urysohn and sequential properties, even for uncountable posets.

The order-theoretic convergence description then interacts directly with the Scott topology. Our main sequential duality theorem states that, for an arbitrary poset $P$, the topologies $\omega(P)$ and $\sigma(P)^d$ have the same convergent sequences, with the same limits. The topologies themselves can differ \cite{Xu-duals}, but they coincide whenever the lower topology is sequential. In particular, topology-level equality follows if the lower topology is Fr\'echet--Urysohn or first-countable, and it holds for every countable poset. Thus the global duality studied in the literature has an unconditional sequential analogue: ordinary sequences cannot distinguish $\omega(P)$ from the de Groot dual of the Scott topology.

For coherent Scott spaces we then determine the lower closure of a $k$-irreducible compact saturated set $K$:
\[
 \cl_{\omega(P)}K=K^{\ell u},
\]
where $K^\ell$ is the set of lower bounds of $K$ and $K^{\ell u}$ is the set of upper bounds of $K^\ell$. If $\wedge K$ exists, the formula reduces to
\[
 \cl_{\omega(P)}K={\uparrow}\!\wedge K.
\]
Together with lower sequential closedness of compact saturated sets, this yields an infimum criterion for co-sobriety. We apply the criterion to bounded-complete dcpos. If the lower topology is sequential---in particular, if it is Fr\'echet--Urysohn or first-countable---then the Scott space is co-sober. Hence the first-countability criterion gives a purely order-theoretic sufficient condition without requiring the dcpo itself to be countable.

On the other hand, Scott co-sobriety can fail without a countability assumption even for complete Boolean algebras. Ji and Xu constructed a complete Boolean algebra whose Scott space is not co-sober \cite{Ji-Xu-Boolean}. Their example is necessarily uncountable, since every countable complete Boolean algebra is finite \cite[Proposition 5.4]{Ji-Xu-Boolean}. Thus the uncountable and countable settings exhibit sharply different behavior. One of the main consequences below is therefore a structural theorem for countable ordered sets with sufficient completeness: every countable bounded-complete dcpo has a co-sober Scott space, and hence so does every countable complete lattice.

Miao, Xi, Li and Zhao constructed a countable complete distributive lattice whose Scott space is not sober \cite{Miao-Xi-Li-Zhao}. Our result shows that their example is nevertheless co-sober. Thus, even among countable complete distributive lattices, Scott sobriety and Scott co-sobriety separate.

The paper is organized as follows. Section 2 develops the order theory of sequential convergence in the lower topology, including the lower-limsup and coinitiality criteria for first countability, and then connects this order convergence with compact saturated subsets of Scott spaces. Section 3 establishes sequential duality between the lower topology and the co-compact topology. Section 4 derives the lower-closure and infimum criteria for coherent Scott spaces. Section 5 applies these results to bounded-complete dcpos and complete lattices.

\section{Sequential convergence in the lower topology}

For $A\subseteq P$, put
\[
 {\uparrow} A=\{x\in P:(\exists a\in A)\ a\le x\},
 ~
 {\downarrow} A=\{x\in P:(\exists a\in A)\ x\le a\}.
\]
For a singleton, we write ${\uparrow} x$ and ${\downarrow} x$.

The \emph{lower topology} $\omega(P)$ is generated by
\[
 \{P\setminus{\uparrow} x:x\in P\}
\]
as a subbase. Thus the finite intersections
\[
 P\setminus{\uparrow} F
 :=P\setminus\bigcup_{a\in F}{\uparrow} a,
 ~ F\subseteq P\text{ finite},
\]
form a base.

A set $U\subseteq P$ is \emph{Scott open} if it is an upper set and, whenever a directed set $D$ has a supremum with $\vee D\in U$, one has $D\cap U\ne\emptyset$. The Scott open sets form the topology $\sigma(P)$, called the \emph{Scott topology} of $P$.

For a $T_0$-space $X$, let $Q(X)$ be the family of its nonempty compact saturated subsets. No separation axiom is included in the meaning of compactness. A member $K\in Q(X)$ is $k$-\emph{irreducible} if
\[
 K=K_1\cup K_2,~ K_1,K_2\in Q(X),
\]
implies $K=K_1$ or $K=K_2$. The space $X$ is \emph{co}-\emph{sober} if every $k$-irreducible member of $Q(X)$ is the saturation of a point. A $T_0$-space is \emph{coherent} if the intersection of any two compact saturated subsets is compact.

A subset $A$ of a topological space $X$ is said to be \emph{sequentially closed} if for any sequence in $A$, $A$ contains all its limits. The space is \emph{sequential} if every sequentially closed subset is closed. It is known that every first-countable space is Fr\'echet--Urysohn and every Fr\'echet--Urysohn space is sequential (see, e.g., \cite[Theorem 1.6.14]{Engelking1989}).

For $A\subseteq P$, define
\[
 A^\ell=\{x\in P:x\le a\text{ for every }a\in A\},
\]
and
\[
 A^{\ell u}=\{y\in P:x\le y\text{ for every }x\in A^\ell\}.
\]

\begin{proposition}[Order characterization of lower convergence]\label{prop:lower-convergence}
Let $(x_n)_{n\ge1}$ be a sequence in a poset $P$ and let $x\in P$. The following are equivalent.
\begin{enumerate}[label=\rm(\arabic*)]
\item $x_n\to x$ in $\omega(P)$.
\item For every $a\in P$ with $a\nleq x$, one has $a\nleq x_n$ eventually.
\item
\[
 \bigcap_{m\ge1}{\downarrow}\{x_n:n\ge m\}\subseteq{\downarrow} x.
\]
\end{enumerate}
\end{proposition}

\begin{proof}
A subbasic lower-open set $P\setminus{\uparrow} a$ contains $x$ exactly when $a\nleq x$. Convergence to all subbasic neighbourhoods is equivalent to convergence in the topology, because eventual membership in finitely many subbasic neighbourhoods implies eventual membership in their intersection. This proves (1) $\Leftrightarrow$ (2).

Assume (2), and let $z$ belong to the intersection in (3). If $z\nleq x$, then (2) gives $z\nleq x_n$ eventually. On the other hand, membership of $z$ in every tail-downset means that for every $m$ there is $n\ge m$ with $z\le x_n$, a contradiction. Thus $z\le x$.

Conversely, if (2) fails, there is $a\nleq x$ such that $a\le x_n$ for infinitely many $n$. Hence for every $l$ some $n\ge l$ satisfies $a\le x_n$, so
\[
 a\in\bigcap_{l\ge1}{\downarrow}\{x_n:n\ge l\},
\]
contradicting (3). Therefore (3) $\Rightarrow$ (2).
\end{proof}

For a sequence $(x_n)$ in a complete lattice $L$, define its \emph{lower limsup} by
\[
 \lsupl x_n
 :=\bigvee\bigcap_{m\ge1}{\downarrow}\{x_n:n\ge m\}.
\]
We also write
\[
 \lsupo x_n:=\bigwedge_{m\ge1}\bigvee_{n\ge m}x_n
\]
for the usual lattice limsup.

\begin{corollary}\label{cor:lower-limsup}
Let $L$ be a complete lattice and $(x_n)_{n\geq 1}$ a sequence in $L$. Then, for every $x\in L$,
\[
 x_n\to x\text{ in }\omega(L)
 ~\Longleftrightarrow~
 \lsupl x_n\le x.
\]
Moreover,
\[
 \lsupl x_n\le \lsupo x_n.
\]
\end{corollary}

\begin{proof}
By Proposition \ref{prop:lower-convergence}, lower convergence to $x$ is equivalent to
\[
 \bigcap_{m\ge1}{\downarrow}\{x_n:n\ge m\}\subseteq{\downarrow} x,
\]
which, in a complete lattice, is equivalent to the join of the left-hand side being at most $x$.

Let
\[
 A=\bigcap_{m\ge1}{\downarrow}\{x_n:n\ge m\}.
\]
For $a\in A$ and every $m$, some $n\ge m$ satisfies $a\le x_n$, hence
\[
 a\le\bigvee_{n\ge m}x_n.
\]
Therefore $a\le\lsupo x_n$ for every $a\in A$. Taking joins gives $\lsupl x_n\le\lsupo x_n$.
\end{proof}

\begin{proposition}\label{prop:chain-limsup}
If $L$ is a complete chain, then for every sequence $(x_n)_{n\geq 1}$ in $L$,
\[
 \lsupl x_n=\lsupo x_n.
\]
Consequently,
\[
 x_n\to x\text{ in }\omega(L)
 ~\Longleftrightarrow~
 \lsupo x_n\le x.
\]
\end{proposition}

\begin{proof}
Put
\[
 A=\bigcap_{m\ge1}{\downarrow}\{x_n:n\ge m\},
 ~
 a=\vee A,
 ~
 r=\bigwedge_{m\ge 1}\bigvee_{n\ge m}x_n.
\]
By Corollary \ref{cor:lower-limsup}, $a\le r$. Suppose that $a<r$.

If there is $b\in L$ with $a<b<r$, then for every $m$ we have
\[
 b<r\le\bigvee_{n\ge m}x_n.
\]
Since $L$ is a chain, some $n\ge m$ must satisfy $b\le x_n$; otherwise $b$ would be an upper bound of the tail. Hence $b\in A$, contradicting $b>a=\vee A$.

It remains to consider the case where no element lies strictly between $a$ and $r$. For every $m$, since $a<r\le\bigvee_{n\ge m}x_n$, some $n\ge m$ satisfies $a<x_n$. By the absence of an element strictly between $a$ and $r$, this implies $r\le x_n$. Thus $r\in A$, again contradicting $a=\vee A<r$. Hence $a=r$.
\end{proof}

\begin{remark}\label{rem:real-limsup}
For the extended real line $\overline{\mathbb R}$ with its usual order, Proposition \ref{prop:chain-limsup} says
\[
 x_n\to x\text{ in the lower topology}
 ~\Longleftrightarrow~
 \limsup_{n\to\infty}x_n\le x.
\]
Thus lower-topological convergence is a direct order-topological extension of the familiar one-sided limsup condition from analysis.
\end{remark}

We next isolate the exact local countability condition for the lower topology. A subset $C$ of a poset $S$ is called \emph{coinitial} if for every $s\in S$ there is $c\in C$ with $c\le s$, or equivalently, $S={\uparrow} C$.

\begin{theorem}[First-countability criterion]\label{thm:first-countable}
Let $P$ be a poset and $x\in P$. The lower topology $\omega(P)$ is first-countable at $x$ if and only if $P\setminus{\downarrow} x$ has a countable coinitial subset.
\end{theorem}

\begin{proof}
Suppose first that $\omega(P)$ is first-countable at $x$. Let $(B_n)_{n\ge1}$ be a countable local base at $x$. For every $n$, choose a basic lower-open neighbourhood
\[
 V_n=P\setminus{\uparrow} F_n,
\]
where $F_n\subseteq P\setminus{\downarrow} x$ is finite, such that
\[
 x\in V_n\subseteq B_n.
\]
Put $C=\bigcup_{n\ge1}F_n$. Then $C$ is countable and $C\subseteq P\setminus{\downarrow} x$.

Let $a\in P\setminus{\downarrow} x$. Since $P\setminus{\uparrow} a$ is a neighbourhood of $x$, there is $n$ with
\[
 B_n\subseteq P\setminus{\uparrow} a.
\]
Hence
\[
 P\setminus{\uparrow} F_n=V_n\subseteq P\setminus{\uparrow} a,
\]
so ${\uparrow} a\subseteq{\uparrow} F_n$. In particular, $a\in{\uparrow} F_n$, and therefore some $c\in F_n\subseteq C$ satisfies $c\le a$. Thus $C$ is coinitial in $P\setminus{\downarrow} x$.

Conversely, let $C\subseteq P\setminus{\downarrow} x$ be countable and coinitial. If $C=\emptyset$, then $P\setminus{\downarrow} x=\emptyset$, so $x$ is a greatest element and $\{P\}$ is a local base at $x$. Assume therefore that $C\ne\emptyset$, and enumerate it as $C=\{c_1,c_2,\ldots\}$, repeating elements if $C$ is finite. Put
\[
 U_n=P\setminus\bigcup_{k=1}^n{\uparrow} c_k.
\]
Then each $U_n$ is a neighbourhood of $x$. We show that $(U_n)_{n\geq 1}$ is a local base at $x$. Let
\[
 V=P\setminus\bigcup_{i=1}^m{\uparrow} a_i
\]
be a basic neighbourhood of $x$, so each $a_i\nleq x$. By coinitiality, choose $c_{k_i}\le a_i$ for $i=1,\ldots,m$. If $n\ge\max_i k_i$, then
\[
 \bigcup_{i=1}^m{\uparrow} a_i
 \subseteq
 \bigcup_{k=1}^n{\uparrow} c_k,
\]
and hence $U_n\subseteq V$. Thus $\omega(P)$ is first-countable at $x$.
\end{proof}

\begin{corollary}\label{cor:first-countable-global}
For a poset $P$, the following are equivalent.
\begin{enumerate}[label=\rm(\arabic*)]
\item $\omega(P)$ is first-countable.
\item For every $x\in P$, the poset $P\setminus{\downarrow} x$ has a countable coinitial subset.
\end{enumerate}
In either case $\omega(P)$ is Fr\'echet--Urysohn and hence sequential.
\end{corollary}

We now turn to compact saturated subsets of Scott spaces. The following compact-down lemma is standard; see, for example, \cite[Section II-1]{Gierz-et-al} and \cite[Lemma 3.14]{Xu-Zhao-Rudin}. The proof is included because the topology on $P$ is otherwise arbitrary.

\begin{lemma}\label{lem:compact-down}
Let $P$ be equipped with a topology in which every principal upper set is closed. If $A\subseteq P$ is compact in that topology, then ${\downarrow} A$ is Scott closed.
\end{lemma}

\begin{proof}
The set ${\downarrow} A$ is a lower set. Let $D\subseteq{\downarrow} A$ be directed, suppose that $d=\vee D$ exists, and assume that $d\notin{\downarrow} A$. Then $d\nleq a$ for every $a\in A$. For each $a\in A$, choose $d_a\in D$ such that $d_a\nleq a$. Thus
\[
 \{P\setminus{\uparrow} e:e\in D\}
\]
is an open cover of $A$. Compactness gives $e_1,\ldots,e_n\in D$ with
\[
 A\subseteq\bigcup_{i=1}^n(P\setminus{\uparrow} e_i).
\]
Choose $e\in D$ above all $e_i$. Since $P\setminus{\uparrow} e_i\subseteq P\setminus{\uparrow} e$ for every $i$, we have $A\subseteq P\setminus{\uparrow} e$. On the other hand, $e\in{\downarrow} A$, so $e\le a$ for some $a\in A$, a contradiction. Hence $d\in{\downarrow} A$.
\end{proof}

\begin{lemma}\label{lem:tails}
Suppose that $x_n\to x$ in $\omega(P)$. For $n\ge1$, let
\[
 A_n=\{x\}\cup\{x_m:m\ge n\},
 ~
 F_n={\downarrow} A_n.
\]
Then every $F_n$ is Scott closed and
\[
 \bigcap_{n\ge1}F_n={\downarrow} x.
\]
\end{lemma}

\begin{proof}
The set $A_n$ is a convergent sequence tail together with its limit, hence it is compact in $\omega(P)$. Principal upper sets are closed in the lower topology, so Lemma \ref{lem:compact-down} shows that $F_n$ is Scott closed.

The inclusion ${\downarrow} x\subseteq\bigcap_nF_n$ is immediate. Conversely, let $z\in\bigcap_nF_n$ and suppose that $z\nleq x$. Then $x\in P\setminus {\uparrow}z\in \omega(P)$. As $x_n\to x$ in $\omega(P)$, there is $m$ such that $x_n\in P\setminus {\uparrow}z$ for all $n\geq m$, and hence $z\not\in F_m$, which contradicts $z\in\bigcap_nF_n\subseteq F_m$. Thus $\bigcap_{n\ge1}F_n={\downarrow} x$.
\end{proof}

\begin{proposition}\label{prop:compact-seq-closed}
Every compact saturated subset of $\Sigma P$ is sequentially closed in $\omega(P)$.
\end{proposition}

\begin{proof}
Let $K\in Q(\Sigma P)$, and suppose that $x_n\in K$ and $x_n\to x$ in $\omega(P)$. Let $F_n$ be as in Lemma \ref{lem:tails}. The sets $K\cap F_n$ form a decreasing sequence of nonempty closed subsets of the compact subspace $K$: indeed, $x_n\in K\cap F_n$. Hence
\[
 \bigcap_{n\ge1}(K\cap F_n)\ne\emptyset.
\]
No separation axiom is used in this compactness argument. By Lemma \ref{lem:tails}, the intersection is $K\cap{\downarrow} x$. Choose $z\in K\cap{\downarrow} x$. Since $K$ is an upper set and $z\le x$, it follows that $x\in K$.
\end{proof}

\section{Sequential duality with the co-compact topology}

Let $\sigma(P)^d$ be the topology generated by
\[
 \{P\setminus K:K\in Q(\Sigma P)\}
\]
as a subbase.

\begin{theorem}[Sequential duality theorem]\label{thm:sequential-duality}
For every poset $P$, the topologies $\omega(P)$ and $\sigma(P)^d$ have the same convergent sequences and the same limits.
\end{theorem}

\begin{proof}
For every $p\in P$, the principal upper set ${\uparrow} p$ is compact and saturated in $\Sigma P$: a Scott open set containing $p$ contains all of ${\uparrow} p$. Therefore $P\setminus{\uparrow} p$ is open in $\sigma(P)^d$, and
\[
 \omega(P)\subseteq\sigma(P)^d.
\]
Thus convergence in $\sigma(P)^d$ implies convergence, to the same limit, in $\omega(P)$.

Conversely, suppose that $x_n\to x$ in $\omega(P)$, and let $K\in Q(\Sigma P)$ with $x\notin K$. If infinitely many $x_n$ belonged to $K$, a subsequence contained in $K$ would still converge to $x$ in $\omega(P)$. Proposition \ref{prop:compact-seq-closed} would then give $x\in K$, a contradiction. Hence $x_n\notin K$ eventually. The sequence is eventually in every subbasic $\sigma(P)^d$-neighbourhood of $x$, and therefore in every finite intersection of such neighbourhoods. Thus $x_n\to x$ in $\sigma(P)^d$.
\end{proof}

\begin{corollary}\label{cor:topological-duality}
If $\omega(P)$ is sequential, then
\[
 \omega(P)=\sigma(P)^d.
\]
In particular, the equality holds whenever $\omega(P)$ is Fr\'echet--Urysohn or first-countable, and it holds for every countable poset $P$.
\end{corollary}

\begin{proof}
By Theorem \ref{thm:sequential-duality}, every $\sigma(P)^d$-open set is sequentially open with respect to $\omega(P)$. If $\omega(P)$ is sequential, it is therefore $\omega(P)$-open. The reverse inclusion always holds. The Fr\'echet--Urysohn and first-countable cases follow from the standard implications. If $P$ is countable, the defining subbase of $\omega(P)$ is countable, so $\omega(P)$ is second-countable and hence sequential.
\end{proof}

Combining Theorem \ref{thm:first-countable} with Corollary \ref{cor:topological-duality} gives a purely order-theoretic version.

\begin{corollary}\label{cor:coinitial-duality}
If, for every $x\in P$, the set $P\setminus{\downarrow} x$ has a countable coinitial subset, then
\[
 \omega(P)=\sigma(P)^d.
\]
\end{corollary}

\begin{remark}\label{rem:global-v-sequential}
Corollary \ref{cor:topological-duality} concerns the co-compact topology of the Scott space. It does not identify this topology, or the lower topology, with the Scott topology of the order dual $P^{\mathrm{op}}$.

The distinction between Theorem \ref{thm:sequential-duality} and topology-level de Groot duality is also worth emphasizing. Global criteria for $\sigma(P)^d=\omega(P)$ involve additional structural hypotheses and are part of the theory of stable compactness and ordered patch spaces; see \cite{Lawson-stably-compact,Xu-duals,Xu-order-topology}. Theorem \ref{thm:sequential-duality} requires no such hypothesis: it asserts an unconditional equality of the sequential convergence structures. Sequentiality of $\omega(P)$ is exactly what is needed here to promote this sequential duality to equality of the topologies.
\end{remark}

\section{Coherent Scott spaces}

\begin{proposition}\label{prop:lower-closure}
Suppose that $\Sigma P$ is coherent. If $K\in Q(\Sigma P)$ is $k$-irreducible, then
\[
 \cl_{\omega(P)}K=K^{\ell u}.
\]
\end{proposition}

\begin{proof}
For every $a\in K^\ell$, one has $K\subseteq{\uparrow} a$. Since ${\uparrow} a$ is closed in $\omega(P)$,
\[
 \cl_{\omega(P)}K
 \subseteq
 \bigcap_{a\in K^\ell}{\uparrow} a
 =K^{\ell u}.
\]

For the reverse inclusion, let $y\notin\cl_{\omega(P)}K$. There are $a_1,\ldots,a_n\in P$ such that
\[
 y\in P\setminus({\uparrow} a_1\cup\cdots\cup{\uparrow} a_n)
 ~\text{and}~
 K\subseteq{\uparrow} a_1\cup\cdots\cup{\uparrow} a_n.
\]
Hence
\[
 K=\bigcup_{i=1}^n(K\cap{\uparrow} a_i).
\]
Discard the empty terms. Each principal upper set is compact saturated, and coherence implies that every remaining $K\cap{\uparrow} a_i$ is compact saturated. Finite unions of compact saturated sets are compact saturated, so binary $k$-irreducibility applies inductively to this finite union. Thus $K=K\cap{\uparrow} a_i$ for some $i$. Hence $a_i\in K^\ell$. Since $a_i\nleq y$, we have $y\notin K^{\ell u}$.
\end{proof}

\begin{corollary}\label{cor:lower-closure-inf}
Under the hypotheses of Proposition \ref{prop:lower-closure}, if $\wedge K$ exists, then
\[
 \cl_{\omega(P)}K={\uparrow}\!\wedge K.
\]
\end{corollary}

\begin{proof}
If $c=\wedge K$, then $K^\ell={\downarrow} c$ and $K^{\ell u}={\uparrow} c$.
\end{proof}

\begin{theorem}\label{thm:infimum-criterion}
Suppose that $\Sigma P$ is coherent and $\omega(P)$ is sequential. Then $\Sigma P$ is co-sober if and only if every $k$-irreducible member of $Q(\Sigma P)$ has an infimum in $P$.
\end{theorem}

\begin{proof}
Suppose first that $\Sigma P$ is co-sober. Its specialization order is the original order on $P$, so the saturation of $x$ is ${\uparrow} x$. Hence every $k$-irreducible $K\in Q(\Sigma P)$ has the form $K={\uparrow} x$, and $\wedge K=x$.

Conversely, let $K\in Q(\Sigma P)$ be $k$-irreducible and suppose that $\wedge K$ exists. By Proposition \ref{prop:compact-seq-closed}, $K$ is sequentially closed in $\omega(P)$. Since $\omega(P)$ is sequential, $K$ is closed. Corollary \ref{cor:lower-closure-inf} now gives
\[
 K=\cl_{\omega(P)}K={\uparrow}\!\wedge K.
\]
Thus every $k$-irreducible compact saturated set is the saturation of a point.
\end{proof}

\section{Bounded-complete dcpos}

We use the convention that a bounded-complete dcpo is a dcpo in which every subset having an upper bound has a supremum, including the empty subset. In particular, $\bot=\vee\emptyset$ exists. This convention is stated explicitly because a weaker usage of the term also occurs in the literature.

Xi and Lawson proved that the Scott space of every complete lattice or bounded-complete dcpo is well-filtered and coherent \cite[Corollary 4.3]{Xi-Lawson}.

\begin{theorem}\label{thm:bc-sequential}
Let $P$ be a bounded-complete dcpo. If $\omega(P)$ is sequential, then $\Sigma P$ is co-sober.
\end{theorem}

\begin{proof}
We verify the infimum condition in Theorem \ref{thm:infimum-criterion}. Let $A\subseteq P$ be nonempty. Its set of lower bounds
\[
 A^\ell=\{x\in P:x\le a\text{ for every }a\in A\}
\]
is nonempty because $\bot\in A^\ell$. It is bounded above by any fixed $a_0\in A$, so
\[
 c=\vee A^\ell
\]
exists. Every $a\in A$ is an upper bound of $A^\ell$, and therefore $c\le a$. Thus $c$ is the greatest lower bound of $A$. In particular, every $k$-irreducible compact saturated subset has an infimum. Coherence follows from the result of Xi and Lawson, so Theorem \ref{thm:infimum-criterion} applies.
\end{proof}

\begin{corollary}\label{cor:bc-frechet}
Let $P$ be a bounded-complete dcpo. If $\omega(P)$ is Fr\'echet--Urysohn, then $\Sigma P$ is co-sober. In particular, this holds if $\omega(P)$ is first-countable.
\end{corollary}

\begin{proof}
Every Fr\'echet--Urysohn space is sequential, and every first-countable space is Fr\'echet--Urysohn. Apply Theorem \ref{thm:bc-sequential}.
\end{proof}

\begin{corollary}\label{cor:bc-coinitial}
Let $P$ be a bounded-complete dcpo. If, for every $x\in P$, the set $P\setminus{\downarrow} x$ has a countable coinitial subset, then $\Sigma P$ is co-sober.
\end{corollary}

\begin{proof}
By Theorem \ref{thm:first-countable}, $\omega(P)$ is first-countable. Apply Corollary \ref{cor:bc-frechet}.
\end{proof}

\begin{corollary}\label{cor:countable-bc}
Every countable bounded-complete dcpo has a co-sober Scott space.
\end{corollary}

\begin{proof}
The lower topology of a countable poset is second-countable and hence sequential. Apply Theorem \ref{thm:bc-sequential}.
\end{proof}

\begin{corollary}\label{cor:countable-lattice}
Every countable complete lattice, and hence every countable complete Heyting algebra, has a co-sober Scott space.
\end{corollary}

The countability hypothesis in Corollary \ref{cor:countable-lattice} is essential for such a general statement. Indeed, Ji and Xu constructed a complete Boolean algebra $B$ whose Scott space $\Sigma B$ is not co-sober \cite{Ji-Xu-Boolean}. Since every countable complete Boolean algebra is finite, this counterexample is necessarily uncountable. Consequently, even within complete Boolean algebras, Scott co-sobriety is automatic in the countable case but may fail once countability is dropped.

\begin{corollary}\label{cor:co-not-sober}
There exists a countable complete distributive lattice whose Scott space is co-sober but not sober.
\end{corollary}

\begin{proof}
Miao, Xi, Li and Zhao constructed a countable complete distributive lattice whose Scott space is not sober \cite{Miao-Xi-Li-Zhao}. Corollary \ref{cor:countable-bc} gives co-sobriety.
\end{proof}

\section{Concluding remarks}

The results above show how order structure, sequential convergence and Scott topology interact in a direct way. Proposition \ref{prop:lower-convergence} shows that lower convergence is controlled by eventual order avoidance; on complete chains it is exactly a limsup condition. Theorem \ref{thm:sequential-duality} shows that ordinary sequences cannot distinguish the lower topology from the co-compact topology of the Scott space, even though the two topologies may be different. First countability and the Fr\'echet--Urysohn property provide natural bridges from this sequential duality to topology-level duality.

For coherent Scott spaces, Proposition \ref{prop:lower-closure} makes the order-theoretic content of lower closure explicit. When the relevant infima exist, co-sobriety is reduced to the passage from lower sequential closedness to lower closedness. This explains why sequential properties of the lower topology lead directly to co-sobriety results for bounded-complete dcpos.

Several questions suggested by this viewpoint remain open. First, it would be useful to obtain intrinsic order-theoretic characterizations of when $\omega(P)$ is Fr\'echet--Urysohn or sequential, extending the first-countability criterion in Theorem \ref{thm:first-countable}. Second, for complete lattices one may ask for structural conditions under which the lower limsup $\lsupl$ agrees with the lattice limsup $\lsupo$ beyond the case of complete chains. Finally, the sequentiality hypothesis in Theorem \ref{thm:bc-sequential} is sufficient but is not claimed to be necessary; weaker closure assumptions ensuring that $k$-irreducible compact saturated subsets are lower closed would lead to further co-sobriety criteria.

\section*{Statements and Declarations}

\paragraph{Funding.}
This research was supported by the National Natural Science Foundation of China (Nos. 12471070, 12071199, 12561083).

\paragraph{Competing interests.}
The authors have no relevant financial or non-financial interests to disclose.

\paragraph{Use of generative AI.}
During the preparation of this manuscript, the authors used ChatGPT (OpenAI) for language editing, presentation suggestions, and manuscript checking. The authors independently verified all mathematical arguments and references, revised the output as necessary, and take full responsibility for the final content.

\paragraph{Data availability.}
No data were used for the research described in this article.

\end{document}